%
\output={\if N\header\headline={\hfill}\fi
\plainoutput\global\let\header=Y}
\magnification\magstep1
\tolerance = 500
\hsize=14.4true cm
\vsize=22.5true cm
\parindent=6true mm\overfullrule=2pt
\newcount\kapnum \kapnum=0
\newcount\parnum \parnum=0
\newcount\procnum \procnum=0
\newcount\nicknum \nicknum=1
\font\ninett=cmtt9

\font\ninebf=cmbx9

\font\sixbf=cmbx6
\font\ninesl=cmsl9

\font\nineit=cmti9

\font\ninerm=cmr9

\font\sixrm=cmr6
\font\ninei=cmmi9
\font\eighti=cmmi8
\font\sixi=cmmi6
\skewchar\ninei='177 \skewchar\eighti='177 \skewchar\sixi='177
\font\ninesy=cmsy9
\font\eightsy=cmsy8
\font\sixsy=cmsy6
\skewchar\ninesy='60 \skewchar\eightsy='60 \skewchar\sixsy='60
\font\titelfont=cmr10 scaled 1440
\font\paragratit=cmbx10 scaled 1200

\font\name=cmcsc10
\font\emph=cmbxti10

\font\tenmsbm=msbm10
\font\sevenmsbm=msbm7
%

%

%
\font\teneufm=eufm10
\font\seveneufm=eufm7
\font\fiveeufm=eufm5
\newfam\eufmfam
\textfont\eufmfam=\teneufm
\scriptfont\eufmfam=\seveneufm
\scriptscriptfont\eufmfam=\fiveeufm

\font\tenmsam=msam10
\font\sevenmsam=msam7
\font\fivemsam=msam5
\newfam\msamfam
\textfont\msamfam=\tenmsam
\scriptfont\msamfam=\sevenmsam
\scriptscriptfont\msamfam=\fivemsam
\font\tenmsbm=msbm10
\font\sevenmsbm=msbm7
\font\fivemsbm=msbm5
\newfam\msbmfam
\textfont\msbmfam=\tenmsbm
\scriptfont\msbmfam=\sevenmsbm
\scriptscriptfont\msbmfam=\fivemsbm
\def\Bbb#1{{\fam\msbmfam\relax#1}}
\def\cz{{\kern0.4pt\Bbb C\kern0.7pt}
}
\def\az{{\kern0.4pt\Bbb A\kern0.7pt}
}
\def\ez{{\kern0.4pt\Bbb E\kern0.7pt}
}
\def\fz{{\kern0.4pt\Bbb F\kern0.3pt}}
\def\gz{{\kern0.4pt\Bbb Z\kern0.7pt}}
\def\hz{{\kern0.4pt\Bbb H\kern0.7pt}
}
\def\kz{{\kern0.4pt\Bbb K\kern0.7pt}
}
\def\nz{{\kern0.4pt\Bbb N\kern0.7pt}
}
\def\oz{{\kern0.4pt\Bbb O\kern0.7pt}
}
\def\rz{{\kern0.4pt\Bbb R\kern0.7pt}
}
\def\sz{{\kern0.4pt\Bbb S\kern0.7pt}
}
\def\pz{{\kern0.4pt\Bbb P\kern0.7pt}
}
\def\qz{{\kern0.4pt\Bbb Q\kern0.7pt}
}
\newskip\ttglue
\def\ninepoint{\def\rm{\fam0\ninerm}%
  \textfont0=\ninerm \scriptfont0=\sixrm \scriptscriptfont0=\fiverm
  \textfont1=\ninei \scriptfont1=\sixi \scriptscriptfont1=\fivei
  \textfont2=\ninesy \scriptfont2=\sixsy \scriptscriptfont2=\fivesy
  \textfont3=\tenex \scriptfont3=\tenex \scriptscriptfont3=\tenex
  \def\it{\fam\itfam\nineit}%
  \textfont\itfam=\nineit
  \def\sl{\fam\slfam\ninesl}%
  \textfont\slfam=\ninesl
  \def\bf{\fam\bffam\ninebf}%
  \textfont\bffam=\ninebf \scriptfont\bffam=\sixbf
   \scriptscriptfont\bffam=\fivebf
  \def\tt{\fam\ttfam\ninett}%
  \textfont\ttfam=\ninett
  \tt \ttglue=.5em plus.25em minus.15em
  \normalbaselineskip=11pt
  \font\name=cmcsc9
  \let\sc=\sevenrm
  \let\big=\ninebig
  \setbox\strutbox=\hbox{\vrule height8pt depth3pt width0pt}%
  \normalbaselines\rm
  \def\sl{\it}}

\headline={\ifodd\pageno\rightheadline\else\leftheadline\fi}
\def\rightheadline{\ninepoint Paragraphen"uberschrift\hfill\folio}
\def\leftheadline{\ninepoint\folio\hfill Chapter"uberschrift}
\let\header=Y
\def\titel#1{\need 9cm \vskip 2truecm
\parnum=0\global\advance \kapnum by 1
{\baselineskip=16pt\lineskip=16pt\rightskip0pt
plus4em\spaceskip.3333em\xspaceskip.5em\pretolerance=10000\noindent
\titelfont Chapter \uppercase\expandafter{\romannumeral\kapnum}.
#1\vskip2true cm}\def\leftheadline{\ninepoint
\folio\hfill Chapter \uppercase\expandafter{\romannumeral\kapnum}.
#1}\let\header=N
}
\def\Titel#1{\need 9cm \vskip 2truecm
\global\advance \kapnum by 1
{\baselineskip=16pt\lineskip=16pt\rightskip0pt
plus4em\spaceskip.3333em\xspaceskip.5em\pretolerance=10000\noindent
\titelfont\uppercase\expandafter{\romannumeral\kapnum}.
#1\vskip2true cm}\def\leftheadline{\ninepoint
\folio\hfill\uppercase\expandafter{\romannumeral\kapnum}.
#1}\let\header=N
}
\def\need#1cm {\par\dimen0=\pagetotal\ifdim\dimen0<\vsize
\global\advance\dimen0by#1 true cm
\ifdim\dimen0>\vsize\vfil\eject\noindent\fi\fi}
\def\neupara#1{\par\penalty-2000
\procnum=0\global\advance\parnum by 1
\vskip1cm\noindent{\paragratit \the\parnum. #1}%
\def\rightheadline{\ninepoint\S\the\parnum.\ #1\hfill \folio}%
\vskip 8mm\noindent}
\def\Proclaim #1 #2\finishproclaim {\bigbreak\noindent
{\bf#1\unskip{}. }{\it#2}\medbreak\noindent}
%
\gdef\proclaim #1 #2 #3\finishproclaim {\bigbreak\noindent%
\global\advance\procnum by 1
{%
{\relax\ifodd \nicknum
\hbox to 0pt{\vrule depth 0pt height0pt width\hsize
   \quad \ninett#3\hss}\else {}\fi}%
\bf\the\parnum.\the\procnum\ #1\unskip{}. }
{\it#2}
\immediate\write\num{\string\def
 \expandafter\string\csname#3\endcsname
 {\the\parnum.\the\procnum}}
\medbreak\noindent}
\newcount\stunde \newcount\minute \newcount\hilfsvar
\def\uhrzeit{
    \stunde=\the\time \divide \stunde by 60
    \minute=\the\time
    \hilfsvar=\stunde \multiply \hilfsvar by 60
    \advance \minute by -\hilfsvar
    \ifnum\the\stunde<10
    \ifnum\the\minute<10
    0\the\stunde:0\the\minute~Uhr
    \else
    0\the\stunde:\the\minute~Uhr
    \fi
    \else
    \ifnum\the\minute<10
    \the\stunde:0\the\minute~Uhr
    \else
    \the\stunde:\the\minute~Uhr
    \fi
    \fi
    }

 \def\calH{{\cal H}}

\def\dim{\mathop{\rm dim}\nolimits}

\def\GL{\mathop{\rm GL}\nolimits}

\def\Hom{\mathop{\rm Hom}\nolimits}

\def\mod{\mathop{\rm mod}\nolimits}

\def\Sp{\mathop{\rm Sp}\nolimits}

\def\Sym{\mathop{\rm Sym}\nolimits}
\def\boxit#1{
  \vbox{\hrule\hbox{\vrule\kern6pt
  \vbox{\kern8pt#1\kern8pt}\kern6pt\vrule}\hrule}}
\def\Boxit#1{
  \vbox{\hrule\hbox{\vrule\kern2pt
  \vbox{\kern2pt#1\kern2pt}\kern2pt\vrule}\hrule}}
\def\rahmen#1{
  $$\boxit{\vbox{\hbox{$ \displaystyle #1 $}}}$$}

\def\smallni{\smallskip\noindent }

\def\Isom{\mathop{\;{\buildrel \sim\over\longrightarrow }\;}}
\def\lo{\longrightarrow}

\def\loma{\longmapsto}

\def\spitz#1{\langle#1\rangle}

\def\pii{\pi {\rm i}}
\def\veps{\varepsilon}
\def\set#1{\bigl\{\,#1\,\bigr\}}

\def\square{\hbox{\hbox to 0pt{$\sqcup$\hss}\hbox{$\sqcap$}}}
\def\qed{\ifmmode\square\else{\unskip\nobreak\hfil
\penalty50\hskip3em\null\nobreak\hfil\square
\parfillskip=0pt\finalhyphendemerits=0\endgraf}\fi}
\def\pn{\the\parnum.\the\procnum}
\def\downmapsto{{\buildrel
        {\vbox{\hbox{\hskip.2pt$\scriptstyle-$}}}
        \over{\raise7pt\vbox{\vskip-4pt\hbox{$\textstyle\downarrow$}}}}}
\input geemen.num
\nopagenumbers
\def\Sym{{\rm Sym}}
\def\tr{{\rm tr}}
\def\Ind{{\rm Ind}}
\def\eps{{\varepsilon}}
\def\grad{{\rm grad}}
\immediate\newwrite\num
\nicknum=0  

\let\header=N
\immediate\newwrite\num\immediate\openout\num=geemen.num
\def\RAND#1{\vskip0pt\hbox to 0mm{\hss\vtop to 0pt{%
\raggedright\ninepoint\parindent=0pt%
\baselineskip=1pt\hsize=2cm #1\vss}}\noindent}
\noindent
\centerline{\titelfont Some Remarks to a Theorem of van Geemen}%
\def\leftheadline{\ninepoint\folio\hfill
 Some Remarks to a Theorem of van Geemen}%
\def\rightheadline{\ninepoint Introduction\hfill \folio}%
\headline={\ifodd\pageno\rightheadline\else\leftheadline\fi}
\vskip 1.5cm
\leftline{\it \hbox to 6cm{Eberhard Freitag\hss}
Riccardo Salvati
Manni  }
  \leftline {\it  \hbox to 6cm{Mathematisches Institut\hss}
Dipartimento di Matematica, }
\leftline {\it  \hbox to 6cm{Im Neuenheimer Feld 288\hss}
Piazzale Aldo Moro, 2}
\leftline {\it  \hbox to 6cm{D69120 Heidelberg\hss}
 I-00185 Roma, Italy. }
\leftline {\tt \hbox to 6cm{freitag@mathi.uni-heidelberg.de\hss}
salvati@mat.uniroma1.it}
\vskip1cm
\centerline{Heidelberg-Roma 2025}
\vskip1cm\noindent
{\paragratit Introduction}%
\vskip1cm\noindent
We consider the Riemann theta series
$$\vartheta[m](\tau,z)=\displaystyle
\sum_{p\in\gz^g}\exp\Bigl(\pii\tau\Bigl[p+{a\over 2}\Bigr]
+2\Bigl(p+{a\over 2}\Bigr)'\Bigl(z+{b\over 2}\Bigr)\Bigr).$$
Here $m\in\gz^{2g}$ is a column vector that is divided into two columns 
$a,b\in\gz^{g}$,
$$m=\pmatrix{a\cr b}.$$
The matrix $\tau$ is from the Siegel upper half plane $\calH_g$. This means that
$\tau$ is a symmetric $g\times g$-matrix such that its imaginary part is positive
definite. Finally $z$ is a column in $\cz^{g}$. Of particular interest are the 
theta nullwerte
$$\vartheta[m](\tau):=\vartheta[m](\tau,0).$$
Up to sign they depend only on $m$ mod two. In this context the columns
$m\in\gz^{2g}$ are called characteristics. A characteristic $m$ is called even
(odd) if $a'b\equiv0$ mod 2 ($\equiv1$ mod 2).
One knows that a theta nullwert vanishes if and only if $m$ is odd.
\smallskip
Let $T_g$ be the vector space generated by the fourth powers of the theta nullwerte.
These are modular forms of weight 2 with respect to the principal congruence
subgroup of level two. 
 In [Ge], Bert van Geemen's result states
\rahmen{\dim (T_g)=(2^g+1)(2^{g-1}+1)/3.}
In  [SM2], it has  been  observed, as  consequence 
of results of Fay [Fa],  that  all linear relations between 
the  $\vartheta_m^4$ are consequences of the quartic Riemann relations.
\smallskip
In this  note, we  want to give a new proof of these results and extend them.
Starting from Fay's results  about  Riemann's relations [Fa],
we  prove that the representation of $\Sp(g,\fz_2)$,  induced by a certain 
character of 
the theta group, has  two  irreducible components  that  
result to be  isomorphic  to the eigenspaces   of Fay's matrix
$M^+$. One  component is  isomorphic to the space 
of Riemann's relations, the other is  isomorphic to $T_g$.
\smallskip
Then we prove a similar result for vector valued modular forms that are attached
to odd characteristics. Here we have to use Fay's matrix $M^-$.
\smallskip
In a last  section  we treat the linear dependencies between arbitrary
powers $\vartheta[m]^k$. We will show that $k=4$ is the only case where
such dependencies can occur.
\smallskip
It is a pleasure to acknowledge the discussions we had with   Bert van Geemen.
 \neupara{The symplectic group over the field of two elements}%
We use the following notations for matrices $A$. 
The transposed of $A$ is denoted
by $A'$.  Its trace is denoted by $\tr(A)$ and $A_0$ is the 
diagonal
of $A$ written
as column vector.
 Let $A$ be an $n\times n$-matrix and $B$ an $n\times m$-matrix, then
$$A[B]:=B'AB.$$
Let $G=G_g=\Sp(g,\fz_2)$ be the symplectic group over the field of two elements.
We recall the affine action of $G$ on $\fz^{2g}$,
$$\sigma\{m\}=\sigma'^{-1}m+\pmatrix{(CD')_0\cr (AB')_0.}$$
Here
$$\sigma=\pmatrix{A&B\cr C&D}$$ 
denotes the decomposition of $\sigma$ into
4 $g\times g$-blocks. It is an action in the sense
$$(\sigma\tau)\{m\}=\sigma\{\tau\{m\}\}.$$
In this context the elements of $\fz_2^{2g}$ are also called
characteristics. A  characteristic
$$m=\pmatrix{a\cr b}$$
is called even if $a'b=0$ otherwise odd. 
We denote the subset of even (odd) characteristics by
$$(\fz_2^{2g})_{{\rm even}},\quad (\fz_2^{2g})_{{\rm odd}}.$$
Their orders are 
$$k_g^{\pm}=  2^{g-1}(2^g\pm1).$$
It is known that the affine action preserves
the parity. Moreover, the group $G$ acts double transitively on the set of
all even (odd) characteristics and it acts transitively on the set of all pairs $(m,n)$
where $m$ is even and $n$ odd (see [Ig1], Chapt.~V, Prop.~2).
\smallskip
We introduce the stabilizers
$$H(m)=\set{\sigma\in G;\quad \sigma\{m\}=m}.$$
In the case $m=0$ this is the so-called theta group. We denote it by $H=H(0)$.
We will have to use the following function on $G\times \fz_2^{2g}$,
$$\eps_m(\sigma)=(-1)^{\tr(B'C)+(B'D)[a]+(A'C)[b]}.$$
This function  arises in the theta transformation formalism as we
will see in the next section. At the moment we take it as it is. 
It has three important properties.
\proclaim
{Lemma}
{The following rules hold.
\item{\rm a)} 
One has that $\eps_m(\sigma)$ for fixed $m$ is a non-trivial (quadratic)
character on $H(m)$. For $h\in H$ one has $\eps_m(h)=\eps_0(h)$.
\item{\rm b)}
For $\sigma,\tau\in G$ one has
$$\eps_m(\tau\sigma)=\eps_m(\sigma)\eps_{\sigma\{m\}}(\tau).$$
\item{\rm c)}
For $h\in H$, $\sigma\in G$ one has
$$\eps_m(h)\eps_m(\sigma)=\eps_m(h\sigma).$$
}
Lvh%
\finishproclaim
{\it Proof.} It can be proved directly by means of somewhat tedious calculations.
Since it is also a direct consequence of the theta transformation formalism which we
treat in the next section, we skip the direct proof and leave it to the interested reader
as an exercise.
\qed
\proclaim
{Lemma}
{The following rule holds.
$$e(\sigma\{m\},\sigma\{n\})=e(m,n)\eps_m(\sigma)e_n(\sigma).$$
}
Lvk%
\finishproclaim
{\it Proof.}
We use the notations
$$\sigma=\pmatrix{A&B\cr C&D},\quad\sigma(m):=\sigma'^{-1}m$$
and
$$m=\pmatrix{a\cr b},\quad n=\pmatrix{\alpha\cr\beta},\quad
v = \pmatrix {(CD')_0\cr  (AB')_0 }.$$
In the proof we will use
$$\eqalign{
(D'B)_0+D'(AB')_0+B'(CD')_0&=0,\cr
(A'C)_0+C'(AB')_0+A'(CD')_0&=0.\cr}$$
This follows from $\sigma^{-1}\{\sigma\{0\}\}=0$.
Now we get
$$\eqalign{
e ( \sigma\{m\},  \sigma\{n\}) 
&= e(\sigma(m), \sigma(n)) e( \sigma(m), v) e( v, \sigma(n))\cr
 &= e(m, n) e(\sigma(m+n), v)\cr
 & = e(m, n)  (-1)^{ (CD')_0'(B(a+\alpha) + 
 A(b+\beta)-   (AB')_0'(D(a+\alpha)+C(b+\beta)}\cr
 &= e(m, n)  (-1)^{ (D'B)_0 ' (a+ \alpha) + (A'C)_0 '( b+\beta)}\cr
 &= e(m, n) \eps_m(\sigma) \eps_n(\sigma)\cr
 }$$
which proves the Lemma.
\qed\smallskip
As special case of Lemma \Lvh\ we get a character 
$$\eps_H:=\eps_0:H\lo \{\pm 1\}.$$
We induce the character $\eps_H$ on $H$ to a $G$- representation 
$\Ind_H^G(\eps_H)$.
So $\Ind_H^G(\eps_H)$ consists of all functions $f:G\to\cz$ with the property
$f(hg)=\eps_H(h)f(g)$. The group $G$ acts by translation from the right on this space.
We mention that this representation has been studied by Frame [Fra]. He proved
that it decomposes into two irreducible non-isomorphic representations and he 
determined the dimensions of the two components. We will reprove this result with
a completely different method. 
\smallskip
We will construct a basis $X(m)$ of  $\Ind_H^G(\eps_H)$. It is parameterized  by even
characteristics. Notice that the number of even characteristics equals the index
of $H$ in $G$ which is the dimension of  $\Ind_H^G(\eps_H)$. So, let $m$ be an even
characteristic. We choose
$\sigma\in G$ such that $\sigma\{m\}=0$. Then we define the following function,
$$X(m):G\lo \cz.$$
It is defined to be 0 outside the coset $H\sigma$, and on this coset it is defined
by
$$X(m)(h\sigma)=\eps_H(h)\eps_m(\sigma).$$
It is essential that this is independent on the choice of $\sigma$. The proof follows
immediately from Lemma \Lvh.
\smallskip
We associate to each even $m$ a variable $e_m$ and consider the vector space
$$\sum_{m\;{{\rm even}}} \cz e_m.$$
We use the basis above to define an isomorphism
$$\Ind_H^G(\eps_H)\Isom \sum_{m\;{{\rm even}}} \cz e_m.$$
Just map $X(m)$ to $e_m$. 
\smallskip
Next we describe an action of $G$ on the right hand side such that this isomorphism
is $G$-equivariant. We have to associate to each $\sigma\in G$ an isomorphism
$$\sigma:\sum_{m\;{{\rm even}}} \cz e_m\Isom \sum_{m\;{{\rm even}}} \cz e_m.$$
On basis elements it is defined by
$$\sigma(e_m)=\eps_m(\sigma)e_{\sigma\{m\}}.$$
It is easy to check that this defines a representation of $G$ and that the isomorphism
$$\Ind_H^G(\eps_H)\Isom \sum_{m\;{{\rm even}}} \cz e_m$$
is $G$-equivariant.
\neupara
{An operator of Fay}%
We treat the space
$$\sum_{m\;{{\rm even}}} \cz e_m$$
together with its $G$-action.
We use the notations
$$\spitz{m,n}=a'\beta+b'\alpha,\qquad m=\pmatrix{a\cr b},\
n=\pmatrix{\alpha\cr \beta}$$
and 
$$e(m,n)=(-1)^{\spitz{m,n}}.$$
Notice that $\spitz{m,n}$ is the symplectic pairing whose invariance group is
$G=\Sp(g,\fz_2)$. Fay [Fa] defined the following operator.
$$M^+: \sum_{m\;{{\rm even}}} \cz e_m\lo  \sum_{m\;{{\rm even}}} \cz e_m,
\quad e_m\loma \sum_n e(m,n)e_n.$$
(Actually Fay ordered the even characteristics and introduced $M^+$ as a matrix,
see [Fa], formula (6).)
\proclaim
{Theorem (Fay)}
{The  operator  $M^+$
has two eigen values, namely $-2^{g-1}$ and  $2^g$. Let
$$V^+\subset\cz^{k_g^+},\quad W^+\subset\cz^{k_g^+}$$
be the eigenspaces. Then we have
$$\cz^{k_g^+}=V^+\oplus W^+.$$
The dimension of the eigenspaces are
$$\quad \dim V^+=(2^{2g}-1)/3,\quad \dim W^+=(2^g+1)(2^{g-1}+1)/3.$$
}
Tfay%
\finishproclaim
We need the following result.
\proclaim
{Lemma}
{The operator $M^+$ commutes with the action of $G$.}
Lmp%
\finishproclaim
{\it Proof.} Use Lemma \Lvk.
\qed\smallskip
As a consequence the two eigenspaces are $G$-invariant. So we see that
$\Ind_H^G(\eps_H)$ admits two invariant subspaces, one isomorphic to $V^+$ the 
other to $W^+$.
\smallskip
We need the decomposition of $G$ into double cosets with respect to $H$.
This can be found also in [Fra].
\proclaim
{Lemma}
{There are two double cosets in $H\backslash G/H$.}
Ldc%
\finishproclaim
We present a different proof.  Let $g_1,g_2\in G-H$. We have to show that they
are in the same double coset. We consider the two characteristics $m_1=g_1\{0\}$,
$m_2=g_2\{0\}$. They are both different from zero. The group $H$ acts transitively
on the set of non-zero characteristics (since $G$ acts double transitively on the set
of even characteristics). Hence there exist $h_1\in H$ such that
$h_1\{m_1\}=m_2$. This implies $h_2:=g_2^{-1}h_1g_1\in H$. It follows
$g_2=h_1g_1h_2^{-1}$. 
\qed
\proclaim
{Lemma}
{The representation 
$\Ind_H^G(\eps_H)$ is either irreducible or decomposes into two 
non-isomorphic irreducible components.}
Lrd%
\finishproclaim
{\it Proof.} The statement is equivalent to
$$\dim\Hom_G(\Ind_H^G(\eps_H),\Ind_H^G(\eps_H))\le2.$$
We use Frobenius reciprocity.
For a representation of $G$ on a finite dimensional vector space $V$
the evaluation map
$$\Hom_G(V,\Ind_H^G(\eps_H))\Isom \Hom_H(V,\cz(\eps_H))$$
is an isomorphism.  Here $\cz(\eps_H)$ is the one dimensional vector space 
$\cz$ where
$H$ acts by means of the character $\eps_H$.
So  we have to show
$$\dim\Hom_H(\Ind_H^G(\eps_H),\cz(\eps_H))\le 2.$$
Since both representations are real, they are isomorphic to their duals.
Hence it  is equivalent to prove
$$\dim\Hom_H(\cz(\eps_H),\Ind_H^G(\eps_H))\le2.$$
Such a homomorphism is defined by 
the image $f$ of $1\in\cz$ and $f$ must have the property
$$f(hg)=\eps_H(h)f(g)\quad\hbox{and}\quad f(gh)=\eps_H(h)f(g).$$
Since $G$ has only two double cosets mod $H$, the dimension of the space
of all $f$ is $\le2$.
\qed\smallskip
We found two non-zero invariant subspaces of $\Ind_H^G(\veps)$. So 
Lemma \Lrd\ shows
that they are irreducible and non-isomorphic.
So we have proved the following theorem.
\proclaim
{Theorem (Frame)}
{The representation $\Ind_H^G(\eps_H)$ decomposes into two non-isomorphic
representations of dimensions
$$(2^{2g}-1)/3,\quad (2^g+1)(2^{g-1}+1)/3.$$
}
TF%
\finishproclaim
\neupara{Theta series}%
The symplectic group $\Sp(g,\rz)$ acts on the Siegel half plane $\calH_g$
through
$$\sigma(\tau)=(A\tau+B)(C\tau+D)^{-1}.$$
We recall   that the Riemann theta functions  are defined by
\smallni
\def\d{\displaystyle}
$$\vartheta[m](\tau,z)=
\d\sum_{p\in\gz^g}\exp\Bigl(\pii\tau\Bigl[p+{a\over 2}\Bigr]
+2\Bigl(p+{a\over 2}\Bigr)'\Bigl(z+{b\over 2}\Bigr)\Bigr).$$
Here
$$m=\pmatrix{a\cr b}\in\gz^{2g}$$
is an integral column vector. We denote by
$$\Gamma_g=\Sp(g,\gz)$$
the Siegel modular group. We consider the map
$$\Gamma_g\times \gz^{2g}\lo \gz^{2g},\quad 
\sigma\{m\}=\sigma'^{-1}m+\pmatrix{(CD')_0\cr (AB')_0.}$$
It is compatible with the map
$\Sp(g,\fz_2)\times\fz_2^{2g}\lo\fz_2^{2g}$ that we defined in Sect.~1.
Hence we can use the same notation.
\proclaim
{Theorem}
{We have
$$\eqalign{
\vartheta[\sigma\{ m\}]&(\sigma( \tau),(C\tau+D)'^{-1}z)=
 \exp\pii\bigl\{((C\tau+D)^{-1}C)[z]\bigr\}\cr
& v(\sigma, m)\sqrt{\det(C\tau+D)}
 \vartheta[ m](\tau,z)\cr}$$
for all $\sigma\in\Gamma_n$. Here
$v(\sigma, m)$ denotes a system of complex numbers which
is independent of $(\tau,z)$.
One  has
$v(\sigma, m)^8=1$.
\smallni
{\bf Corollary.} For the theta nullwerte this means
$$
\vartheta[\sigma\{ m\}](\sigma(\tau))=
 v(\sigma, m)\sqrt{\det(C\tau+D)}
 \vartheta[ m](\tau).$$
}
Ttt%
\finishproclaim
The system of numbers $v(\sigma,m)$ is very complicated. We will need only the
4th powers. It is known that these depend only on $m$ mod 2 and on the image
of $\sigma$ with respect to the homomorphism $\Sp(g,\gz)\to\Sp(g,\fz_2)$.
This homomorphism is surjective. So we can define
$$v(\sigma,m)^4\quad\hbox{for}\quad \sigma\in\Sp(g,\fz_2),\ m\in\fz_2^{2g}.$$
This function is known explicitly.
\proclaim
{Lemma}
{One has 
$$v(\sigma,m)^4=\eps_m(\sigma)=(-1)^{\tr(B'C)+(B'D)[a]+(A'C)[b]}.$$
}
Lveps%
\finishproclaim
A proof can be found in [Ig2], Lemma 10.
\qed\smallskip
This was the expression we introduced in Sect.~1. It gets now a natural interpretation.
Lemma \Lvh\  is an  easy consequences of this interpretation.
\smallskip
For $\sigma\in\Gamma_g$ and a function $f:\calH_g\to\cz $ and an 
integer $r$ we will use the notation.
$$(f\vert_r\sigma)(\tau)=\det(C\tau+D)^{-r}(f(\sigma\tau).$$
The theta transformation formula shows that
$$f\mapsto f\vert_2 \sigma^{-1}$$ 
defines a representation of $\Gamma_g$ on
$T_g$. This representation factors through $G=\Sp(g,\fz_2)$.
Now we consider the linear map
$$\Ind_H^G(\eps_H)\lo T_g,\quad X(m)\loma\vartheta[m]^4.$$
The theta transformation formula shows that this map is compatible with the action
of $G$. It is also clear that it is surjective. Its dimension is smaller than that  of
$\Ind_H^G(\eps_H)$.  To show this, one needs that there exists at least one non-trivial
relation between the $\vartheta[m]^4$. The quartic Riemann relations give one.
(In the
 next Theorem
we will describe a space of such relations.)
Therefore one of the two irreducible components must map
isomorphically to $T_g$. We have to work out which one.
\smallskip
For this we formulate some of the Riemann theta relations in a way that 
has been described by
Fay [Fa]. In Sect.~1  we introduced  the operator
$$M^+: \sum_{m\;{{\rm even}}} \cz e_m\lo  \sum_{m\;{{\rm even}}} \cz e_m$$
and described the eigenspaces. The original meaning of this operator is also
due to Fay [Fa]. He showed that the Riemann quartic relations imply the
following result.
\proclaim
{Theorem (Fay)}
{Let $v$ be an eigen vector of $M^+$ with eigenvalue $-2^{g-1}$ (i.e.~$v\in V^+$).
Then the relation
$$\sum_m v(m)\vartheta[m]^4=0$$
holds.
}
Tfayz%
\finishproclaim
This shows that the irreducible component that maps isomorphically 
to $V^+$ goes
to zero in $T_g$. Hence the other component $W^+$ maps isomorphically to
$T_g$.
We obtain the first main result.
\proclaim
{Theorem}
{We have
$$\dim T_g=(2^g+1)(2^{g-1}+1)/3.$$
The linear relations between the $\vartheta[m]^4$ are consequences of
Riemann's quartic relations.}
TvG%
\finishproclaim
\neupara{A  variant of van Geemen's result}%
Let $f(\tau,z)$ be a holomorphic function on $\calH_g\times \cz^g$.
Then we can define the holomorphic vector valued function
$$J(f):\calH_g\lo\cz^g,\quad J(f)=\grad_z f(\tau,z)\vert_{z=0}.$$
We are interested in $J(\vartheta[m])$. This is different from 0 if and only if
$m$ is odd. 
\proclaim
{Theorem}
{For $\sigma\in \Gamma_g$ and odd $m\in\gz^{2g}$ the transformation formula
$$J(\vartheta[\sigma\{m\}])(\sigma(\tau))=
v(M,\sigma)\sqrt{\det(C\tau+D)}(C\tau+D)(\vartheta[m])(\tau)$$
holds.}
Tot%
\finishproclaim
This is a consequence of Theorem \Ttt\ (see for example [SM1]).
\qed\smallskip
We consider the fourth symmetric tensor product  $\Sym^4(\cz^g)$ and define the 
functions
$$S(\vartheta[m]):\calH_g\to \Sym^4(\cz^g),\quad 
\tau\loma \Sym^4J(\vartheta[m])(\tau).$$
The space generated by these $k_g^-$ tensors is called  by $S_g$.
A special case of Theorem \Tot\ states
$$S(\vartheta[\sigma\{m\}])(\sigma(\tau))=
\eps_m(\sigma){\det(C\tau+D)}^{2}\Sym^4(C\tau+D)S(\vartheta[m])(\tau).$$
This shows that the $S(\vartheta[m])$ are vector valued modular forms
on the principal congruence subgroup of level two with the character $\eps_0$
and with respect to the $\GL(n,\cz)$ representation
$$\det(CZ+D)^2\Sym^4(CZ+D).$$
This representation is irreducible and has highest weight 
$(6,2,\dots,2)$.
\smallskip
Up to the weight this is the same transformation formula as for the
$\vartheta[m]^4$ in the even case.
\smallskip
The tensors $S(\vartheta[m])$ satisfy similar linear relations 
as the $\vartheta[m]^4$.
Instead of the theta group $H$ we now have to use the stabilizer $K=H(n)$ of 
a fixed odd characteristic $n$. To be concrete, we take
$$n'=(1,0,\dots,0,1,0,\dots,0).$$
Now we have to consider the non trivial character $\eps_K=\eps_n$. The induced
representation $\Ind_K^G(\eps_K)$ comes into the game. Its dimension
equals the number $k_g^-$ of odd characteristics.
We will construct a basis $Y(m)$ of  $\Ind_K^G(\eps_K)$. 
It is parameterized  by odd
characteristics. So, let $m$ be an odd
characteristic. We choose
$\sigma\in G$ such that $\sigma\{m\}=n$. Then we define the following function,
$Y(m):G\lo \cz$ analoguos to $X(m)$.
It is defined to be 0 outside the coset $K\sigma$, and on this coset it is defined
by
$$X(m)(k\sigma)=\eps_K(k))\eps_m(\sigma).$$
It is also independent of the choice of $\sigma$. Again this follows
from Lemma \Lvh.
\smallskip
We associate to each odd $m$ a variable $e_m$ and consider the vector space
$$\sum_{m\;{{\rm odd}}} \cz e_m.$$
We use the basis above to define an isomorphism
$$\Ind_K^G(\eps_K)\Isom \sum_{m\;{{\rm odd}}} \cz e_m.$$
Just map $Y(m)$ to $e_m$. 
\smallskip
Besides $M^+$ Fay [Fa] defined the following operator
$$M^- \sum_{m\;{{\rm odd}}} \cz e_m\lo  \sum_{m\;{{\rm odd}}} \cz e_m,
\quad e_m\loma \sum_n e(m,n)e_n.$$
\proclaim
{Theorem (Fay)}
{The  operator  $M^-$
has two eigen values, namely $2^{g-1}$ and  $-2^g$. Let
$$V^-\subset\cz^{k_g^+},\quad W^-\subset\cz^{k_g^+}$$
be the eigenspaces. Then we have
$$\cz^{k_g^-}=V^-\oplus W^-.$$
The dimension of the eigenspaces are
$$\quad \dim V^-=(2^{2g}-1)/3,\quad \dim W^-=(2^g-1)(2^{g-1}-1)/3.$$
}
Tfay%
\finishproclaim
The same proof as in Lemma \Ldc\ shows that $M^-$ commutes with $G$. Hence
the two spaces $V^-$ and $W^-$ are $G$-invariant subspaces. 
\smallskip
With the  same method, using Fay's matrix $M^-$ instead of $M^+$, one can
get now a new proof of the following result of Frame.
\proclaim
{Theorem (Frame)}
{Let $K=H(m)$ be the stabilizer of an odd characteristic $m$. Denote by
$\eps_K$   the 
distinguished character of $K$. Then the representation 
$\Ind_K^G(\veps_K)$ decomposes into two non-isomorphic irreducibles
representation of degrees
$$(2^{2g}-1)/3,\quad (2^g-1)(2^{g-1}-1)/3.$$
}
Rfr%
\finishproclaim
As we have seen, the group $\Sp(g,\fz_2)$ acts one the space $S_g$. The linear
map
$$\Ind_K^G(\eps_K)\lo S_g,\quad Y(m)\loma S(\vartheta[m]),$$
is $G$-linear. It is surjective. Hence it is an isomorphism of one of the two
irreducible complements maps isomorphically. We will show that the
second case holds. For this we have to  make use of Fay's theta relations. A special
case of [Fa]  (23) says for $z_1-z_2=z_3-z_4$ and $w\in W^-$
$$\sum_{m\;{\rm odd}} w_m
\vartheta[m](\tau,z_1)\vartheta[m](\tau,z_2)
\vartheta[m](\tau,z_3)\vartheta[m](\tau,z_4)=0$$
or
$$\sum_{m\;{\rm odd}} w_m
\vartheta[m](\tau,z_1)\vartheta[m](\tau,z_2)
\vartheta[m](\tau,z_3)\vartheta[m](\tau,-z_1+z_2+z_3)=0$$
with three independent variables $z_1,z_2,z_3$. We take the gradient with respect
to the variable $z_3$ and take then $z_3=0$. Since $\vartheta[m](\tau,0)=0$
for odd $m$, we obtain
$$\sum_{m\;{\rm odd}} w_m
\vartheta[m](\tau,z_1)\vartheta[m](\tau,z_2)
\vartheta[m](\tau,-z_1+z_2)J(\vartheta[m])(\tau)=0.$$
Taking the $z_2$-gradient gives
$$\sum_{m\;{\rm odd}} w_m
\vartheta[m](\tau,z_1)
(\vartheta[m](\tau,-z_1)J(\vartheta[m](\tau))\otimes
J(\vartheta[m])(\tau)=0.$$
The theta functions $\vartheta[m](\tau,z)$ are invariant under $z\mapsto-z$
up to a sign. Now differentiating twice by $z_1$ and then setting $z_1=0$ gives
$$\sum_{m\;{\rm odd}} w_m S(\vartheta[m])=0.$$
So we get non trivial linear relations between the $S(\vartheta[m])$. 
\proclaim
{Theorem}
{The dimension of the space spanned by the tensors $S(\vartheta[m])$ has 
dimension $(2^{2g}-1)/3$.}
Smt%
\finishproclaim
\neupara{Arbitrary powers of the theta nullwerte}%
The question arises whether one can get similar results for the space
$T_g(k)$  generated 
by $k$th powers of theta nullwerte where $k$ is a natural number (so $T_g=T_g(4)$).
We start with the case $k\equiv 4\mod 8$. 
In this case we have a natural surjective $G$-linear map
$$\Ind_H^G(\eps_H)\lo T_g(k)$$
that generalizes that described in Sect.~3.
\proclaim
{Theorem}
{Assume $k\equiv 4\mod 8$, $k\ne 4$.
The map
$$\Ind_H^G(\eps_H)\lo T_g(k)$$
is an isomorphism.
Hence $T_g(k)$ is the  direct sum of the two irreducible representations
of dimensions
$$(2^{2g}-1)/3,\quad (2^g+1)(2^{g-1}+1)/3.$$}
Tevk%
\finishproclaim
{\it Proof.} We have to show that none of the two irreducible components goes to 0.
Actually we find in each of the two irreducible components one element that
does not go to 0. We use the description of $V^+,W^+$ as in Theorem \Tfay.
Clearly
$$\eqalign{
(2^g-1)X(0)^k-\sum_{m\ne 0}X(m)^k&\in V^+,\cr
(2^g+1)X(0)^k+\sum_{m\ne 0}X(m)^k&\in W^+.\cr
}$$
We want to show that their images are not 0. This means
$$\eqalign{
(2^g-1)\vartheta[0]^k-\sum_{m\ne 0}\vartheta[m]^k&\ne 0,\cr
(2^g+1)\vartheta[0]^k+\sum_{m\ne 0}\vartheta[m]^k&\ne0.\cr
}$$
We apply to the two functions the Siegel $\Phi$-operator $g-1$ times.
In the first case one gets
$$2^{g-1}\biggl(\vartheta\left[\matrix{0\cr 0}\right]^k-
(\vartheta\left[\matrix{0\cr 1}\right]^k-\vartheta\left[\matrix{1\cr 0}\right]^k
\biggr)
$$
This is different from zero for $k\ne 4$ since
$$\vartheta\left[\matrix{0\cr 0}\right]^4=
\vartheta\left[\matrix{0\cr 1}\right]^4+\vartheta\left[\matrix{1\cr 0}\right]^4
$$
is the defining relation between the three theta nullwerte in genus 1.
The second case is similar.
\qed\smallskip
Next we consider the  case $k\equiv0\mod 8$. Here we have an obvious surjective
$G$-invariant map
$$\Ind_H^G(1)\lo T_g(k)$$
where $1$ denotes the trivial character of $H$. The same argument as in the case
of the nontrivial character $\veps$ shows that this representation has
at most two irreducible components. One is the trivial one dimensional 
representation consisting of all constant functions $f:H\to\cz$. So we have
two irreducible components, one of dimension $1$, the other of dimension
$k_g^+-1$. The representation of dimension $k_g^+-1$ cannot map to zero,
since the dimension of $T_g(k)$ is greater than one. The function 
$\sum_m\vartheta[m]^k$ is non-zero and is invariant under $G$. 
Hence we obtain the following result.
\proclaim
{Theorem}
{Assume $k\equiv 0\mod 8$. Then the map
$$\Ind_H^G(1) \lo T_g(k)$$
is an isomorphism. As a consequence $T_g(k)$ decomposes 
under the action of $G$ into two irreducible
components of dimension
$$1,\quad 2^{g-1}(2^{g}+1)-1.$$
}
Rkze%
\finishproclaim
Finally we treat the case that $k$ is not divisible by $4$. Since we want to include
also odd $k$ we have to replace $\Gamma_g$ by its two fold metaplectic
cover
$\tilde\Gamma_g$. Recall that this is the set of all pairs
$(\sigma,\sqrt{\det(C\tau+D)})$, $\sigma\in\Gamma_g$,
 where $\sqrt{\det(C\tau+D)}$ is one of the 
holomorphic square roots of $\det(C\tau+D)$. 
We set for a function $f:\calH_g\to \cz$ 
$$f\vert(\sigma,\sqrt{\det(C\tau+D)})=
\sqrt{\det(C\tau+D)}\,\strut^{-k}f(\sigma\tau).$$
Then $f\mapsto f\vert\sigma^{-1}$ gives a representation of $\tilde\Gamma_g$
on $T_g(k)$.
If $k$ is even, it factors through $\Gamma_g$ and through $G$ if $k$ is divisible
by 4.
We denote by $\tilde\Gamma_g[2]$ the inverse image of $\Gamma_g[2]$
in $\tilde\Gamma_g[2]$. The group $\tilde\Gamma_g$ and hence 
$\tilde\Gamma_g[2]$ act on the space $T_g(k)$. 
\proclaim
{Lemma}
{Assume that $4$ does not divide $k$.  The space $T_g(k)$ decomposes
under $\tilde\Gamma_g[2]$ into the one dimensional spaces
$\cz\vartheta[m]^k$ which are pairwise non-isomorphic (their characters are
pairewise different).
}
Lnv%
\finishproclaim
{\it Proof.}  We know that  the $\vartheta[m]^k$ are modular forms 
on $\tilde\Gamma_g[2]$  with respect to some character $v_m$. 
We have to show that
the characters $v_m$, $v_n$ for $m\ne n$ are different. Since $\Gamma_g$ acts
doubly transitive, we can assume $m=0$ and $n={a\choose 0}$,  $a\ne 0$. 
One checks that $v_m$ and $v_n$ disagree on translations matrices
${E \,2S\choose 0\; E}$.
\qed
\proclaim
{Theorem}
{Assume that $4$ does not divide $k$. Then the space $T_g(k)$ is irreducible
under the action of $\tilde\Gamma_g$. Its dimension is
$$k_g^+=2^{g-1}(2^{g}+1).$$}
Tvnd%
\finishproclaim
{\it Proof.} Let $V\subset T_g(k)$ be a non-zero invariant subspace. 
Under the action
of $\tilde\Gamma_g[2]$ it decomposes into certain sums of the
one-dimensional spaces
$\cz\vartheta[m]$. The group $\tilde\Gamma$ acts transitively  on the 
even characteristics $m$. Hence $V=T_g(k)$.
\qed\smallskip
An immediate consequence of Theorems \Tevk, \Rkze, \Tvnd\ is the following
result.
\proclaim
{Corollary}
{The theta nullwerte $\vartheta[m]^k$ are linear independent with one exception,
$k=4$.}
Ci%
\finishproclaim
\vskip1cm\noindent
{\paragratit References}%
\bigskip
\item{[Fa]} Fay, J.F.:
{\it On the Riemann-Jacobi Formula,}
Nachr. Akad.  Wiss. G\"ottingen, 1979
\medskip
\item{[Fr]} Frame, J.S.: {\it Some Characters of Orthogonal Groups over the field of
two Elements,}  Lecture Notes in Mathematics {\bf 372}, Proceedings of the Second
International Conference on the Theory of groups, Edited by  M.F.\ Newman
Springer-Verlag Berlin Heidelberg u.a.,\ 1974
\medskip
\item{[Ge]} van Geemen, B.:
{\it Siegel modular forms vanishing on the moduli space of curves,}
Invent. math. {\bf 78}, 1984
\medskip
\item{[Ig1]} Igusa, J.I.:
{\it Theta Functions,}
Grundlehren der mathematischen Wis\-sen\-schaf\-ten, Band {\bf 194},
Springer-Verlag Berlin Heidelberg u.a.,\ 1972
\medskip
\item{[Ig2]} Igusa, J.I.:
{\it On Jacobi's Derivative Formula and its Generalizations,}
Am. Journ. Math., Vol, {\bf 102}, No. 2,  1980
\medskip
\item{[SM1]} Salvati Manni, R.: 
{\it On the dimension of the vector space $\cz[\theta_m]_4$,}
Nagoya Math. J. Vol. {\bf 98}, 1985
\medskip
\item{[SM2]} Salvati Manni, R.: 
{\it On the Not Identically Zero Nullwerte of Jacobians of Theta Functions
with Odd Characteristics,}
Advances in Math. {\bf 47}, 1983
\bye